\documentclass[a4paper]{scrartcl} 

\usepackage[latin1]{inputenc}
\usepackage[ngerman, english]{babel}

\usepackage{amssymb}
\usepackage{amsmath}
\usepackage{amsthm}
\usepackage{mathtools}

\title{The Pólya sum process: Limit theorems for conditioned random fields}
\author{Mathias Rafler\thanks{rafler@ma.tum.de}\\
\small{TU M\"unchen}\\[-7pt]
\small{Zentrum Mathematik, M5}\\[-7pt]
\small{Boltzmannstr. 3}\\[-7pt]
\small{D-85747 Garching bei M\"unchen}
}

\renewcommand{\Pr}{\mathbb P} 
\newcommand{\Poy}{\mathsf S\!} 
\newcommand{\PoyS}[2][\ast]{\mathsf S^{#1}_{\! #2}} 
\newcommand{\Pp}{\mathsf P} 
\newcommand{\B}{\mathcal{B}} 
\newcommand{\Bbd}{\mathcal{B}_{0}} 
\newcommand{\Hsig}{\mathcal{H}}
\newcommand{\Gsig}{\mathcal{G}}
\newcommand{\Fsig}{\mathcal{F}}
\newcommand{\Esig}{\mathcal{E}}

\newcommand{\Htail}{\Hsig_{\infty}}
\newcommand{\Gtail}{\Gsig_{\infty}}

\newcommand{\Etail}{\Esig_{\infty}}

\newcommand{\Hf}{\mathbb H}
\newcommand{\Gf}{\mathbb G}
\newcommand{\Ff}{\mathbb F}
\newcommand{\Ef}{\mathbb E}

\newcommand{\M}{\mathcal M}
\newcommand{\MX}{\mathcal M(X)}

\newcommand{\Mpm}{\mathcal M^{\cdot\cdot}}
\newcommand{\MpmX}{\mathcal M^{\cdot\cdot}(X)}
\newcommand{\MN}{\mathcal M(\N)}
\newcommand{\MpmN}{\mathcal M^{\cdot\cdot}(\N)}
\newcommand{\N}{\mathbb N}
\newcommand{\R}{\mathbb R}
\newcommand{\T}{\mathcal T}
\renewcommand{\d}{\mathrm{d}}
\renewcommand{\exp}{\operatorname{exp}}
\newcommand{\e}{\operatorname{e}}
\newcommand{\id}{\operatorname{id}}
\renewcommand{\phi}{\varphi}
\renewcommand{\theta}{\vartheta}
\newcommand{\eps}{{\mathbf\varepsilon}}
\newcommand{\F}{\mathcal{F}}
\newcommand{\fsm}{\text{-a.s.}}
\newcommand{\supp}{\operatorname{supp}}
\newcommand{\ra}{\rightarrow}
\theoremstyle{definition}
\newtheorem{defdefinition}{Definition}
\theoremstyle{plain}
\newtheorem{defsatz}[defdefinition]{Theorem}

\newtheorem{defprop}[defdefinition]{Proposition}
\newtheorem{deflemma}[defdefinition]{Lemma}
\newtheorem{defmlemma}[defdefinition]{Main Lemma}
\newtheorem{deffolgerung}[defdefinition]{Corollary}
\theoremstyle{remark}
\newtheorem{defbemerkung}[defdefinition]{Remark}

\newcommand{\satzn}[2]{\begin{defsatz}[#1]#2\end{defsatz}}
\newcommand{\prop}[1]{\begin{defprop}#1\end{defprop}}
\newcommand{\propn}[2]{\begin{defprop}[#1]#2\end{defprop}}

\newcommand{\lemma}[1]{\begin{deflemma}#1\end{deflemma}}
\newcommand{\mlemma}[1]{\begin{defmlemma}#1\end{defmlemma}}
\newcommand{\lemman}[2]{\begin{deflemma}[#1]#2\end{deflemma}}
\newcommand{\korollar}[1]{\begin{deffolgerung}#1\end{deffolgerung}}

\numberwithin{equation}{section}

\begin{document}

\maketitle

\begin{abstract}
In \cite{hZ09}, Zessin constructed the so-called Pólya sum process via partial
integration technique. This process shares some important properties with the
Poisson process such as complete randomness and infinite divisibility. This work
discusses H-sufficient statistics for the Pólya sum process as
it was done for the Poisson process in~\cite{hZ76}.
\end{abstract}

\section{Introduction}

Given the number of points to be placed in a bounded region of a polish space, the Poisson process places these points independently and identically distributed. In replacing this
mechanism by a Pólya urn-like scheme one obtains the Pólya sum process. This point process was constructed in~\cite{hZ09} as the unique solution of the integral equation involving the Campbell measure 
\begin{equation} \label{eq:intro:polya}
  \int h(x,\mu) C_{\Pp}(\d x,\d\mu) = z\iint h(x,\mu+\delta_{x})\bigl(\rho+\mu\bigr)(\d x)\Pp(\d\mu).
\end{equation}
With $\rho$ instead of $z(\rho+\mu)$ this is Mecke's characterisation of the Poisson process with intensity measure $\rho$. The additional summand imports a reinforcement.

For the Poisson process with intensity measure $\rho$ this means that independent of a realized point configuration, a new point is added with intensity $\rho$. In contrast, the mechanism of the Pólya sum process rewards points contained in a configuration with their multiplicity, hence forces clustering of points. A similar effect show bosonic particles, and recently in~\cite{BZ11}, the distribution of the particles of the ideal Bose gas on the possible states is connected with the Pólya sum process. This point of view on the ideal Bose gas is different from its position distribution derived in ~\cite{kF80}.

In~\cite{hZ76}, Nguyen and Zessin characterised the mixed Poisson processes as canonical Gibbs states, and among those the Poisson process as the extremal ones. In the following we address a similar question for the Pólya sum process for three different local specifications yielding different families of mixed Pólya sum processes. We identify the set of extreme points as two one- and one two-parameter family of Pólya sum processes.

This result paves the way to a Bayesian viewpoint on the mixed Pólya sum processes along the lines of~\cite{GW82}: Since the Gibbs states are mixed Pólya sum processes directed by some probability measure on a particular parameter set, one might estimate the distribution of these parameters. Indeed, as will turn out, the distribution of the parameters given a single observation, the posteriori measure is concentrated on a single point. This is the characterization of ergodically decomposable priori measures obtained Gl\"otzl and Wakolbinger. 

We follow the approach to construct the Martin-Dynkin boundary as in Dynkin~\cite{eD71a,eD71b} and Foellmer~\cite{hF75}. Terms and notations are adopted to those in Dynkin~\cite{eD78}.

In section~\ref{sect:results} the basic setup is given including the Pólya sum process and its representation in~\ref{subsect:results:psp}. Thereafter
in~\ref{subsect:results:Hsuff} the definition of local specification and H-sufficient statistics is recalled. Finally in~\ref{subsect:results:bd} we construct the local
specifications and give the main results, which are the Martin-Dynkin boundaries of the three local specifications. Their proofs are contained in section~\ref{sect:proofs}.


\section{Pólya sum process and results\label{sect:results}}

\subsection{Pólya sum process\label{subsect:results:psp}}
Let $X$ be a polish space and denote by $\B=\B(X)$ its Borel sets as well as by $\Bbd=\Bbd(X)$ the ring of bounded Borel sets of $X$. Furthermore let $\MX$ and $\MpmX$ be the space of locally finite measures and locally finite point measures on $X$, respectively, each of which is vaguely polish, the $\sigma$-algebras $\F$ generated by the evaluation mappings $\zeta_{B}(\mu)=\mu(B)$, $B\in\B_{0}$. For $B\in\Bbd$ denote by $\F_B$ the $\sigma$-algebra generated by $\zeta_{B'}$ for all $B'\in\Bbd$ such that $\B'\subseteq B$, i.e. the $\sigma$-algebra of inside-$B$ events. We call a probability measure $\Pp$ on $\MX$ a random measure and if $\Pp$ is concentrated on $\MpmX$ a point process. Finally let $F(X)$ be the set of bounded, non-negative and measurable functions on $X$ and $F_b(X)\subset F(X)$ the subset of those functions in $F(X)$ with bounded support.

By $C_{\Pp}$ we denote the Campbell measure of $\Pp$
\begin{equation*}
  C_{\Pp}(h)=\iint_{X\times\MX}h(x,\mu)\mu(\d x)\Pp(\d\mu).
\end{equation*}
$C_{\Pp}$ determines $\Pp$ uniquely on $\MX\setminus\{0\}$ and therefore on $\MX$ since $\Pp$ is a law, see e.g. \cite{MKM78} for details corresponding the Campbell measure. Of particular interest have been disintegrations of the Campbell measure of the type of equation~\eqref{eq:results:Polya-int} below, e.g. in~\cite{MWM79}, \cite{oK78} and~\cite{NZ10}. Recently, Zessin gave a construction method for the reverse direction in~\cite{hZ09}: Given a kernel $\eta$ from $\MpmX$ to $\MX$ is there a point process with Papangelou kernel $\eta$, i.e. satisfies an integral equation of the type~\eqref{eq:results:Polya-int} with $z(\rho+\mu)$ replaced by $\eta$? A particular example he gave is the Pólya sum process $\Poy_{z,\rho}$ for $\rho\in\MX$ and $z\in(0,1)$ which is the unique solution of 
the integral equation
\begin{equation} \label{eq:results:Polya-int}
  C_{\Pp}(h)=z\iint h(x,\mu+\delta_x)\bigl(\rho+\mu\bigr)(\d x)\Pp(\d\mu), \qquad h\in F.
\end{equation}
Directly from equation~\eqref{eq:results:Polya-int} he showed that $\Poy_{z,\rho}$ firstly has independent increments and secondly that the number $\zeta_B$ of points inside the bounded set $B$ obeys a negative binomial distribution. Also, only in using relation~\eqref{eq:results:Polya-int}, one shows that $\Poy_{z,\rho}$ also satisfies
\begin{equation} \label{eq:results:infdiv-int}
  C_{\Pp}(h) = \iint h(x,\mu+\nu)C_L(\d x,\d\nu)\Pp(\d\mu),
\end{equation}
where $L$ is the image of the measure $\rho\otimes\tau_z$ on $X\times\N$ under the mapping $(x,r)\mapsto r\delta_x$ and $\tau_z=\sum_{j\geq 1}\frac{z^j}{j}\delta_j$. This shows that $\Poy_{z,\rho}$ is infinitely divisible with Levy measure $L$. Again, the integral equation~\eqref{eq:results:infdiv-int} determines $\Poy_{z,\rho}$ uniquely.

Since the intensity measure of $\Poy_{z,\rho}$ is
\begin{equation} \label{eq:results:intensity}
  \nu_{\Poy_{z,\rho}}^{1}=\frac{z}{1-z}\rho,
\end{equation}
directly from equation~\eqref{eq:results:infdiv-int} we get that the Palm kernel $\PoyS[x]{z,\rho}$ is given by, where $\Delta_{x}=\delta_{\delta_x}$ is the point process which realizes a point at $x$,
\begin{equation} \label{eq:results:palm-kernel}
  \PoyS[x]{z,\rho} = \Poy_{z,\rho}\ast\frac{1-z}{z}\sum_{j\geq 1}z^j\Delta_x^{\ast j},
\end{equation} 
i.e. the Palm kernel at $x$ is the original process with an additional point with geometrically distributed multiplicity at $x$. A second direct consequence of the integral
equation~\eqref{eq:results:infdiv-int} or~\eqref{eq:results:Polya-int} is that the Laplace functional of $\Poy_{z,\rho}$ is given by
\begin{equation}
  L_{\Poy_{z,\rho}}(f)=\exp\left\{ -\int \log\frac{1-z\e^{-f(x)}}{1-z}\rho(\d x)\right\}.
\end{equation}

During the construction of the Pólya sum process in~\cite{hZ09}, the image of the iterated 
Pólya sum kernels for $B\in\Bbd$, $\rho_B$ denoting the restriction of $\rho$ to
$B$,
\begin{equation}
  \pi^{(m)}_B(0,\d x_1,\ldots \d
  x_m)=\bigl(\rho_B+\delta_{x_1}+\ldots+\delta_{x_{m-1}}\bigr)(\d
  x_m)\cdots\bigl(\rho_B+\delta_{x_1}\bigr)(\d x_2)\rho_B(\d x_1)
\end{equation}
under the mapping $(x_1,\ldots x_m)\mapsto\delta_{x_1}+\ldots+\delta_{x_m}$ occurred. We
keep terms of measures and write
\begin{equation*}
  \Mpm_{(m)}(\N)=\biggl\{ \gamma\in\MpmN:\sum_{j\geq 1}j\gamma(j)=m \biggr\},
\end{equation*}
to denote the set of partitions of $m$ elements represented by the number of
families $\gamma(j)$ of size $j$. For given $\gamma$ denote by the finite product
\begin{equation*}
  c(\gamma)=\prod_{j\in\N}j^{\gamma(j)}\gamma(j)!,
\end{equation*}
then $m!/c(\gamma)$ is the number of permutations with cycle structure given by $\gamma$.
\mlemma{\label{thm:results:iteratedkernel}
Let $\phi$ be non-negative and $\F_B$-measurable. Then
\begin{align*}
\lefteqn{    \int \phi(\delta_{x_1}+\ldots+\delta_{x_{m}}) \pi^{(m)}_B(0,\d x_1,\ldots,\d
    x_m) }\\ 
    &= \sum_{k=1}^{m}\sum_{\stackrel{\gamma\in\Mpm_{(m)}(\N)}{\gamma(\N)=k}}
      \frac{m!}{c(\gamma)} \int_{B^k} \phi(i_1\delta_{x_1}+\ldots +i_k\delta_{x_k})
      \rho(\d x_1)\cdots\rho(\d x_k)\\
    &= \sum_{k=1}^{m}\frac{m!}{k!}\sum_{\stackrel{i_1',\ldots,i_k'\geq 1}{i_1'+\ldots +i_k'=m}}
      \frac{1}{i_1'\cdots i_k'} \int_{B^k} \phi(i_1'\delta_{x_1}+\ldots +i_k'\delta_{x_k})
      \rho(\d x_1)\cdots\rho(\d x_k)      
\end{align*}
where for each fixed $\gamma$ with $\gamma(\N)=k$, $i_1\leq \ldots\leq i_k$ is the unique increasing sequence with exactly $\gamma(j)$ of the $i_k$'s equal to $j$.
}%
One recognises that if $\phi=1$, then the integral on the lhs equals $\rho(B)^{[m]}$
and the inner sum on the rhs is apart from the weight $\rho(B)$ the number of permutations with exacly $k$
cycles ${m \brack k}$, one recovers 
\begin{equation*}
  \sum_{k=1}^{m} \rho(B)^{k}
  \sum_{\stackrel{\gamma\in\Mpm_{(m)}(\N)}{\gamma(\N)=k}} \frac{m!}{c(\gamma)}
  =\rho(B)^{[m]}.
\end{equation*}
The direct consequence is the fact that if $z\in(0,1)$, then
\begin{equation*}
  \begin{multlined}
    \sum_{m\geq 1} \frac{z^m}{m!}\pi^{(m)}_B(0,\phi)\\
    =\sum_{k\geq 1}\frac{1}{k!}\sum_{i_1,\ldots,i_k\geq 1}
      \frac{z^{i_1}\cdot\ldots\cdot z^{i_k}}{i_1\cdot\ldots\cdot i_k}
      \int_{B^k}\phi(i_1\delta_{x_1}+\ldots +i_k\delta_{x_k}) \rho(\d x_1)\cdots\rho(\d x_k).
  \end{multlined}
\end{equation*}

\subsection{H-sufficient statistics and local specifications\label{subsect:results:Hsuff}}

Equipped with $\subseteq$, $\Bbd$ is a partially ordered set. We assume that there is an increasing sequence $(B_n)_{n\in\N}$ of bounded sets such that $\bigcup_{n\in\N} B_n=X$ and for each $B\in\Bbd$ there is a $n_0\in\N$ with $B\subseteq B_{n_0}$. Furthermore we construct a decreasing family $\Ef$ of $\sigma$-fields $\Esig_B\subseteq\Fsig$ indexed by the bounded sets. A family of Markovian kernels $\pi^{\Ef}=(\pi^{\Ef}_{B})_{B\in\Bbd}$ from $\MpmX$ to $\MpmX$ is called a local specification if
\begin{enumerate}
  \item if $B'\subseteq B$, then $\pi^{\Ef}_{B}\pi^{\Ef}_{B'}=\pi^{\Ef}_{B}$; \label{enum:results:locspec:konsistenz}
  \item if $f$ is $\F$-measurable, then $\pi^{\Ef}_{B}f$ is $\Esig_B$-measurable; \label{enum:results:locspec:messbarkeit}
  \item if $f$ is $\Esig_B$-measurable, then $\pi^{\Ef}_{B}f=f$. \label{enum:results:locspec:ident}
\end{enumerate}
Given the family $\pi^\Ef$, one is interested firstly in the convex set $C(\pi^{\Ef})$ of point processes $\Pr$ with the property 
\begin{equation} \label{eq:results:invariant}
  \Pr(\phi|\Esig_B)=\pi^{\Ef}_{B}(\,\cdot\,,\phi)\qquad\Pr\fsm
\end{equation}
and secondly in its extremal points. The former are the $\pi^{\Ef}_{B}$-invariant measures. Since $B\subseteq B_{n_0}$ for each $B\in\Bbd$ and some $n_0$ depending on $B$, $C(\pi^{\Ef})$ agrees with the set of 
all point processes $\Pr$ such that equation~\eqref{eq:results:invariant} holds for each $B_n$.

A $\sigma$-field $\Esig$ is sufficient for  $C(\pi^{\Ef})$ if the $\Pr\in C(\pi^{\Ef})$ have a common conditional distribution given $\Esig$, i.e. there exists $Q_\mu$ such that
\begin{equation*}
  Q_\mu=\Pr\bigl(\,\cdot\,|\Esig\bigr)(\mu)\qquad \Pr\text{-a.e. }\mu.
\end{equation*}
According to~\cite{eD78}, the tail-$\sigma$-field $\Etail=\bigcap_n\Esig_n$ is a H-sufficient statistic for $C(\pi^{\Ef})$, i.e. it is a sufficient statistic and $\Pr\bigl(\mu:Q_\mu\in C(\pi^{\Ef})\bigr)=1$ for all $\Pr\in C(\pi^{\Ef})$. In the situation of this work the family $\pi^{\Ef}=(\pi^{\Ef}_n)_{n}$ given by
\begin{equation*}
  \pi^{\Ef}_n(\mu,\phi)=\Poy_{z,\rho}(\,\cdot\,|\Esig_{B_{n}})(\mu)
\end{equation*}
is a local specification and therefore $\Etail$ is H-sufficient for the set $C(\pi^{\Ef})$ of $\pi^{\Ef}$-invariant point processes. Its essential part $\Delta^{\Ef}$ is the set of extremal points of $C(\pi^{\Ef})$.

Thus the family $\pi^{\Ef}$ describes local laws and the aim is to determine global laws consistent with this description. In particular we get integral representations of elements in $C(\pi^{\Ef})$ in terms of the extremal points $\Delta^{\Ef}$, since the latter set turns out to be a set of Pólya sum processes, we get a characterization of mixed Pólya sum processes.

\subsection{The tail-$\sigma$-fields and results\label{subsect:results:bd}}

Since the multiplicity of points is rather the rule then the exception, there are several possibilities to choose statistics. Having in mind a picture of building
bricks, one may measure on the one hand just the number of sites where they are placed, or on the next hand the total number of bricks which are placed without taking into account
the number of sites, or on the third hand one may measure both. We call these three ensembles the \emph{occupied sites ensemble}, the \emph{total height ensemble} and the
\emph{size-and-height ensemble}. 

We start with the family $\hat{\Ff}=(\hat{\Fsig}_B)_{B\in\Bbd}$ of outside events,
\begin{equation*}
  \hat{\Fsig}_B=\sigma(\zeta_{B,B'}:B\subseteq B'\in\Bbd)
\end{equation*}
for $B\in\Bbd$, where $\zeta_{B,B'}=\zeta_{B'}-\zeta_{B}$ is the increment.

Of first interest is the family $\Gf=(\Gsig_B)_{B\in\Bbd}$ of $\sigma$-fields where $\Gsig_{B}$ is generated by $\hat{\Fsig}_{B}$ and $\sigma(\xi_{B})$, where $\xi_{B}\mu\coloneqq\zeta_{B}\mu^{\ast}$ counts the number occupied sites of the configuration $\mu\in\MpmX$ inside $B\in\Bbd$.

\satzn{Martin-Dynkin boundary of the occupied sites ensemble}{ \label{thm:results:t-spec}
Let $z\in(0,1)$ and $\rho\in\MX$ be a diffuse and infinite measure. The tail-$\sigma$-field $\Gtail$ is H-sufficient for the family 
\begin{equation*}
  C(\pi^{\Gf})=\left\{ \int\Poy_{z,w\rho}V(\d w)| V \text{ probability measure on } [0,\infty) \right\}, 
\end{equation*}
and the set of its extremal points is exactly the family 
\begin{equation*}
  \Delta^{\Gf} = \{ \Poy_{z,w\rho} | 0\leq w<\infty \}.
\end{equation*}
}
Therefore when estimating the number of occupied sites of the particles, we get a one-parameter family. If we replace the number of occupied sites in $B$ by the number of particles in $B$, i.e. if $\Hf=(\Hsig_{B})_{B\in\Bbd}$ is the collection of $\sigma$-algebras $\Hsig_B$ generated by $\hat{\Fsig}_{B}$ and $\sigma(\zeta_{B})$, we obtain

\satzn{Martin-Dynkin boundary of the total height ensemble}{ \label{thm:results:b-spec}
Let $z\in(0,1)$ and $\rho\in\MX$ be a diffuse and infinite measure. The tail-$\sigma$-field $\Htail$ is H-sufficient for the family 
\begin{equation*}
  C(\pi^{\Hf})=\left\{\int \Poy_{z,\rho}V(\d z)| V \text{ probability measure on } [0,1)\right\}
\end{equation*}
and the essential part of the Martin-Dynkin boundary is exactly the family 
\begin{equation*}
  \Delta^{\Hf} = \{ \Poy_{z,\rho} | 0\leq z<1 \}.
\end{equation*}
}
In this case, by estimating the number of particles per volume, we adjust the parameter $z$. Because of~\eqref{eq:results:palm-kernel}, this increases the average multiplicity of the points as well as by~\eqref{eq:results:intensity} the average number of occupied sites. Each of the tail-$\sigma$-fields $\Gtail$ and $\Htail$ is a H-sufficient statistic for an one-parameter family of Pólya sum processes. Finally we combine both statistics and obtain a two-parameter family. Let 
\begin{equation*}
  \Esig_{B} \coloneqq \hat{\Fsig}_{B} \vee \sigma(\xi_{B}) \vee \sigma(\zeta_{B}).
\end{equation*}

\satzn{Martin-Dynkin boundary of the size-and-height ensemble}{ \label{thm:results:tb-spec}
Let $z\in(0,1)$ and $\rho\in\MX$ be a diffuse and infinite measure. The tail-field $\Etail$ is H-sufficient for the family 
\begin{equation*}
  C(\pi^{\Ef})=\left\{\int\Poy_{z,w\rho}V(\d z,\d w)| V \text{ probability measure on } (0,1)\times(0,\infty)\cup\{(0,0)\}\right\},
\end{equation*}
and its extremal points are given exactly by all Pólya sum processes for the pairs $(z,w\rho)$,
\begin{equation*} 
  \Delta^{\Ef}=\{\Poy_{z,w\rho}| 0< z< 1, 0< w<\infty\}\cup\{\delta_0\}.
\end{equation*}
}

Remark that if any of the parameters $z$ and $\rho$ vanishes, then the Pólya sum process realizes the empty configuration almost surely and in this case set both parameters zero.


\section{Proofs \label{sect:proofs}}

In the very first part of this proof section we show the representation of the iterated Pólya sum kernels from lemma~\ref{thm:results:iteratedkernel}. We then turn to the theorems
from section~\ref{subsect:results:bd}.

The basic structure of the proofs of theorems~\ref{thm:results:t-spec}~--~\ref{thm:results:tb-spec} is similar and therefore we start with the general part recalling the lines
of~\cite{eD78,hF75} in section~\ref{subsect:proofs:frame}. The basic problem is to identify the limits $Q$. In section~\ref{subsect:proofs:ose}
theorem~\ref{thm:results:t-spec} is proved by direct computation, theorems~\ref{thm:results:b-spec}~and~\ref{thm:results:tb-spec} are proven by means of large deviations 
in~\ref{subsect:proofs:tpe}~and~\ref{subsect:proofs:tpose}. Their common part, which consists mainly of the proofs of propositions~\ref{thm:proofs:b-lim}~and%
~\ref{thm:proofs:tb-lim} can be found in section~\ref{subsect:proofs:ldp}.

\subsection{The iterated Polya sum kernel}

For any $\gamma\in\Mpm_{(m)}(\N)$ with $\gamma(\N)=k$ define a measure $M^\gamma$ on
$\MpmX$ by
\begin{equation*}
  M^{\gamma}(\phi) = \int \phi(i_1\delta_{x_1}+\ldots +i_k\delta_{x_k}) 
                 \rho(\d x_1)\cdots\rho(\d x_k),\qquad \phi\geq 0,
\end{equation*}
where $i_1\leq \ldots\leq i_k$ is the unique increasing sequence with exactly
$\gamma(j)$ of the $i_k$'s equal to $j$. Then immediatly we get for $\Fsig_B$-measurable, non-negative $\phi$ that
\begin{equation*}
  \int \phi(\delta_{x_1}+\ldots+\delta_{x_{m}}) \pi^{(m)}_B(0,\d x_1,\ldots \d
    x_m)=\sum_{k=1}^{m}\sum_{\stackrel{\gamma\in\Mpm_{(m)}(\N)}{\gamma(\N)=k}} 
    \alpha^{(m)}_\gamma M^\gamma(\phi)
\end{equation*}
for some constants $\alpha^{(m)}_\gamma$ which we identify right after the following lemma
\lemma{
The family of constants $(\alpha^{(m)}_\gamma)_{m\geq 1,\gamma\in\Mpm_{(m)}(\N)}$ satisfies the recursion
\begin{equation} \label{eq:proofs:ipsk:recursion}
  \alpha^{(m+1)}_{\gamma} = \sum_{j\geq 1} j\bigl(\gamma(j)+1\bigr) 
    \alpha^{(m)}_{\gamma+\delta_j-\delta_{j+1}}1_{\{\gamma(j+1)\geq 1\}}
   + \alpha^{(m)}_{\gamma-\delta_{1}} 1_{\{\gamma(1)\geq 1\}}
\end{equation} 
with initial value $\alpha^{(1)}_{\delta_1}=1$.
}
The recursion~\eqref{eq:proofs:ipsk:recursion} states that for given family composition $\gamma\in\Mpm_{(m+1)}(\N)$, $\alpha^{(m+1)}_{\gamma}$ is determined by all possibilities
to introduce a new member to a population of size $m$ weighted with the size of families in the smaller population.

\begin{proof}
Let $m\geq 1$ and $\phi$ be $\F_B$-measurable and non-negative. $\pi^{(m+1)}_B$
puts another point to the realisation of $\pi^{(m)}_B$ either introducing a new
one or putting it to an already existing point. The first case leads to the first summand in the pre-last line with an additional family of size 1. In the second case, putting the new point to an existing one means to add it to an existing family. There are exactly $j\gamma(j)$ members in families of size $j$, hence exactly this number of possibilities to add the new point to a family of size $j$.
\begin{align*}
   \lefteqn{ \int \phi(\delta_{x_1}+\ldots+\delta_{x_{m+1}}) \pi^{(m+1)}_B(0,\d x_1,\ldots,\d x_{m+1})}\\
    &= \iint \phi(\delta_{x_1}+\ldots+\delta_{x_{m+1}}) \bigl(\rho_B+\delta_{x_1}+\ldots
       +\delta_{x_{m}}\bigr)(\d x_{m+1})\pi^{(m)}_B(0,\d x_1,\ldots,\d x_{m})\\
    &=   \begin{multlined}[t]
\sum_{k=1}^{m}\sum_{\stackrel{\gamma\in\Mpm_{(m)}(\N)}{\gamma(\N)=k}} 
       \alpha^{(m)}_{\gamma} \int_{B^{k+1}}\phi(\delta_{x_1}+\ldots+\delta_{x_{m+1}})
       M^{\gamma}(\d x_1,\ldots,\d x_m)\rho(\d x_{m+1})\\
       + \sum_{l=1}^{m}\sum_{k=1}^{m}\sum_{\stackrel{\gamma\in\Mpm_{(m)}(\N)}{\gamma(\N)=k}} 
       \alpha^{(m)}_{\gamma} \int_{B^k} \phi(\delta_{x_1}+\ldots+\delta_{x_{m}}+\delta_{x_l})
       M^{\gamma}(\d x_1,\ldots,\d x_m)
       \end{multlined}\\
    &= \sum_{k=1}^{m}\sum_{\stackrel{\gamma\in\Mpm_{(m)}(\N)}{\gamma(\N)=k}} 
       \alpha^{(m)}_{\gamma}M^{\gamma+\delta_{1}}(\phi)
       + \sum_{k=1}^{m}\sum_{\stackrel{\gamma\in\Mpm_{(m)}(\N)}{\gamma(\N)=k}}\sum_{j\geq 1}
       j\gamma(j)\alpha^{(m)}_{\gamma} M^{\gamma-\delta_j+\delta_{j+1}}(\phi)\\
    &= \sum_{k=1}^{m+1}\sum_{\stackrel{\gamma\in\Mpm_{(m+1)}(\N)}{\gamma(\N)=k}}
       \alpha^{(m+1)}_{\gamma}M^{\gamma}(\phi),
\end{align*}
where the coefficients satisfy the recursion~\eqref{eq:proofs:ipsk:recursion}.
\end{proof}%
Since $\alpha^{(1)}_{\delta_1}=1$, one checks that
$\alpha^{(m)}_{\gamma}=\frac{m!}{c(\gamma)}$ satisfies the
recursion~\eqref{eq:proofs:ipsk:recursion}, which is the first equation in
Main Lemma~\ref{thm:results:iteratedkernel}. The second is an immediate consequence
by noting that
\begin{align*}
   \lefteqn{ \sum_{k=1}^m \sum_{\stackrel{\gamma\in\Mpm_{(m)}(\N)}{\gamma(\N)=k}}
    \frac{m!}{c(\gamma)} M^\gamma(\phi)}\\
     &=
       \sum_{k=1}^m\frac{m!}{k!} \sum_{\stackrel{\gamma\in\Mpm_{(m)}(\N)}{\gamma(\N)=k}}
       \binom{k}{\gamma}\frac{1}{i_1\cdots i_k}
       \int_{B^k} \phi(i_1\delta_{x_1}+\ldots i_k\delta_{x_k})\rho(\d x_1)\cdots\rho(\d x_k)\\
     &=\sum_{k\geq 1}\frac{m!}{k!} \sum_{\stackrel{i_1',\ldots,i_k'\geq
     1}{i_1'+\ldots +i_k'=m}}
       \int_{B^k} \phi(i_1'\delta_{x_1}+\ldots i_k'\delta_{x_k})\rho(\d x_1)\cdots\rho(\d x_k)
\end{align*}
where in the second line the $i_j$'s are ordered and given by $\gamma$, which drops in the last
line since $\phi$ is symmetric when changing the order of that summation.

\subsection{The main frame\label{subsect:proofs:frame}}

We follow the programme given in~\cite{hF75} and give the basic common structure for the local specifications related to $\Gf$, $\Hf$ and $\Ef$. Here we use $\Gf$ and $\pi^\Gf$ to represent any of the three setups. The family $\pi^{\Gf}=(\pi^{\Gf}_B)_{B\in\Bbd}$ with $\pi^{\Gf}_{B}=\Poy_{z,\rho}(\,\cdot\,|\Gsig_{B})$ is a family of Markovian kernels since for $A\in\hat{\Fsig}$, $\Poy_{z,\rho}(A|\Gsig_B)$ is $\sigma(\zeta_B)$-measurable and furthermore a local specification since $\pi^{\Gf}$ clearly satisfies~\ref{enum:results:locspec:konsistenz}~--~\ref{enum:results:locspec:ident}. We want to determine the set
$C(\pi^{\Gf})$ of point processes $\Pr$ which are locally given by
\begin{equation*}
  \Pr(A|\Gsig_B)=\pi^{\Gf}_{B}(\,\cdot\,,A)\qquad \Pr\fsm
\end{equation*}
for all $B\in\Bbd$ and $A\in\Fsig$. $\Poy_{z,\rho}\in C(\pi^{\Gf})$ in each ensemble ensures the non-emptiness of $C(\pi^{\Gf})$. Following Föllmer and Dynkin, the Martin-Dynkin
boundary is constructed in the following way: If $C_{\infty}(\pi^{\Gf})$ is the set of all limits 
\begin{equation*}
  \lim_{k\to\infty}\pi^{\Gf}_{B_k}(\mu_k,\,\cdot\,),
\end{equation*}
where $(\mu_k)_{k\in\N}$ is a sequence in $\Mpm$, then $C_{\infty}(\pi^{\Gf})$ is complete in the set of all probability measures on $\Mpm$, therefore Polish. The measurable
space $C_{\infty}(\pi^{\Gf})$ equipped with the Borel-$\sigma$-field is the Martin-Dynkin boundary of $\pi^{\Gf}$. Since $(B_k)_{k\in\N}$ is an increasing sequence, we have for
each $\Pr\in C(\pi^{\Gf})$ and $\Pr$-integrable $\phi$
\begin{equation*}
  \Pr(\phi|\Gtail)=\lim_{k\to\infty}\pi^{\Gf}_{B_k}(\,\cdot\,,\phi)\qquad \Pr\fsm
\end{equation*}
Therefore we firstly have to compute the $\Pr\fsm$ existing weak limit 
\begin{equation*}
  Q^{\Gf}_{\mu} = \lim_{k\to\infty}\pi^{\Gf}_{B_k}(\mu,\,\cdot\,),
\end{equation*}
which is contained in $C_{\infty}(\pi^{\Gf})$ for $\Pr$-a.e. $\mu$ by construction and in $C(\pi^{\Gf})$ for $\Pr$-a.e. $\mu$ by the H-sufficiency. We will see that $Q^\Gf_\mu$ is a Pólya sum process for $\Pr$-a.e. $\mu$, which implies that
\begin{equation} \label{eq:gen:cond-exp-tail}
  \Pr(\phi|\Gtail)=\Poy_{Z,W\rho}(\phi)\qquad \Pr\fsm
\end{equation}
or
\begin{equation*}
  \Pr(\phi)=\int\Poy_{Z(\mu),W(\mu)\rho}(\phi)\Pr(\d\mu)
\end{equation*}
for suitable, possibly a.s.~constant random variables $Z$ and $W$ on $(\MpmX,\Fsig)$ (even $\Gtail$-measurable) and that $C(\pi^{\Gf})$ consists of mixed Pólya sum processes. Finally we identify the
extremal points $\Delta^{\Gf}$ of $C(\pi^{\Gf})$ as the Pólya sum processes among the mixed ones. 

The important step is to determine the limits $Q^\Gf_\mu$. In the first ensemble $Q_\mu$ can be identified in showing the convergence of Laplace functionals, in the other two we
use a large deviation principle and a conditional minimisation procedure.


\subsection{Occupied sites ensemble \label{subsect:proofs:ose}}

In this first case we are only interested in the number of sites which are occupied, the multiplicity does not matter. The infinite divisibility of $\Poy_{z,\rho}$ admits direct computations. For $B\in\Bbd$, recall that $\xi_B\mu=\zeta_B\mu^\ast$ is the number of points of the support of $\mu$. Then the kernel $\pi^{\Gf}_{B}$ is given by
\begin{equation*}
  \pi^{\Gf}_{B}(\mu,\phi) = \Poy_{z,\rho}\bigl( \phi|\Gsig_{B} \bigr)(\mu) = \Poy_{z,\rho}\bigl( \phi(\,\cdot\, +\mu_{B^{c}}) | \xi_{B}=\xi_{B}\mu\bigr).
\end{equation*}
On $\{\xi_B=n\}$, for $\Fsig_B$-measurable $\phi$ this is because of the diffuseness of $\rho$
\begin{equation} \label{eq:proofs:t-spec}
  \pi^{\Gf}_{B}(\,\cdot\,,\phi) = \begin{multlined}[t]
      \left(\frac{1}{-\log(1-z)\rho(B)}\right)^{n} \\ \times\int_{B^{n}} \sum_{i_{1},\ldots,i_{n}\geq 1} 
      \phi(i_{1}\delta_{x_{1}}+\ldots+i_{n}\delta_{x_{n}}) \frac{z^{i_{1}+\ldots +i_{n}}}{i_1\cdots i_n} \rho(\d x_{1}) \cdots \rho(\d x_{n}).
      \end{multlined}
\end{equation} 
Inside $B$, $\pi^{\Gf}_{B}$ places points at exactly $\xi_B$ sites with independent multiplicities.
We denote by $W_{k}$ the number of occupied sites in $B_{k}$ normalised by $-\log(1-z)\rho(B_{k})$ 
\begin{equation*}
  W_{k}\mu \coloneqq \frac{\xi_{B_{k}}\mu}{-\log(1-z)\rho(B_{k})}.
\end{equation*}
Let 
\begin{equation*}
  \Mpm_\rho=\left\{\mu\in\MpmX: \lim_{k\to\infty} \frac{\xi_{B_{k}}\mu}{-\log(1-z)\rho(B_{k})}\text{ exists}\right\}
\end{equation*}
and set $W$ on $\Mpm_\rho$ the corresponding limit and $0$ otherwise.

\prop{ \label{thm:proofs:t-lim}
Let $\rho\in\MX$ be diffuse and infinite and $z\in(0,1)$. For any $\Pr\in C(\pi^\Gf)$, $\Pr(\M_\rho)=1$ and furthermore for $\phi\in L^{1}(\Pr)$
\begin{equation} \label{eq:proofs:t-lim}
  \Pr\bigl(\phi|\Gtail\bigr)(\mu)=\Poy_{z, W(\mu)\cdot\rho}(\phi) \qquad \Pr\text{-a.e. }\mu.
\end{equation}
}

\begin{proof}
Let
\begin{equation*}
  \Mpm_{\Gf}=\left\{ \mu\in\Mpm: \lim_{k\to\infty}\pi^{\Gf}_{B_k}(\mu,\,\cdot\,) \text{ exists}\right\}.
\end{equation*}
On $\{\xi_B=n\}$ we have for $f\in F_b(X)$, and for $B\in\Bbd$ large enough such that $\supp f\subseteq B$,
\begin{align*}
  \pi^{\Gf}_{B}\Bigl(\,\cdot\,,\e^{-\zeta_f}\Bigr)
    &=\begin{multlined}[t]
       \bigl[-\log(1-z)\rho(B)\bigr]^{-n} \\ \times \int_{B^{n}} \sum_{i_1,\ldots,i_n\geq 1} \e^{-i_1f(x_1)-\ldots-i_nf(x_n)}\frac{z^{i_1+\ldots+i_n}}{i_1\cdots i_n}\rho(\d x_1)\cdots\rho(\d x_n)
       \end{multlined}\\
    &= \left[\frac{1}{-\log(1-z)\rho(B)} \int_B \log(1-z\e^{-f(x)})\rho(\d x)\right]^{n}\\
    &= \left[1-\frac{1}{-\log(1-z)\rho(B)} \int \log\frac{1-z\e^{-f(x)}}{1-z}\rho(\d x)\right]^{\log(1-z)\rho(B)\frac{n}{\log(1-z)\rho(B)}}.
\end{align*}
The restiction of the integral to $B$ can be dropped in the last line since $f(x)\neq 0$ iff $\log\frac{1-z\e^{-f(x)}}{1-z}\neq 0$. If we replace $B$ by $B_k$ such that $\supp f\subseteq B_k$, and let $k\to\infty$, the lhs converges on $\Mpm_{\Gf}$, as well as the inner part of the rhs does to an exponential, and we get that
$\xi_{B_k}\mu/\bigl(-\log(1-z)\rho(B_k)\bigr)$ converges, the limit being $W(\mu)$. Hence $\Mpm_{\Gf}\subseteq\Mpm_{\rho}$. On the contrary, if the limit $W(\mu)$ exists, then the rhs converges and therefore $\Mpm_\rho\subseteq\Mpm_{\Gf}$. Therefore we get
\begin{align*}
  \lim_{k\to\infty}\pi^{\Gf}_{B_k}\Bigl(\,\cdot\,,\e^{-\zeta_f}\Bigr) 
     = \exp\left\{ -W \int_B \log\frac{1-z\e^{-f(x)}}{1-z}\rho(\d x)
     \right\}\qquad \Pr\text{-a.s},
\end{align*}
which is the Laplace functional of the mixed Pólya sum process $\Poy_{z,W\rho}$.
\end{proof}

From~\eqref{eq:proofs:t-lim} we get by taking expectations for bounded and $\Etail$ measurable $\phi$ and bounded and measurable $f:\R_+\to\R$,
\begin{equation*}
  \Pr\bigl(\phi f(W)\bigr) = \Pr\Bigl(\phi\Poy_{z,W\rho}\bigl(f(W)\bigr)\Bigr),
\end{equation*}
hence in particular
\begin{equation*}
  \Poy_{z,W(\mu)}\bigl(W=W(\mu)\bigr)\qquad \Pr\text{-a.e. } \mu.
\end{equation*}
Finally let $V_\Pr$ be the distribution of $W$ under $\Pr$, then by equation~\eqref{eq:proofs:t-lim},
\begin{equation*}
  \Pr(\phi)=\Pr\bigl(\Poy_{z,W\rho}(\phi)\bigr)=\int\Poy_{z,w\rho}(\phi)V_\Pr(\d w)
\end{equation*}
and $\Pr$ is a mixed Pólya sum process. On the other hand, for every probability measure $V$ on $\R_+$ the corresponding mixed Pólya sum process is contained in $C(\pi^\Gf)$. Note that for $\Pr=\Poy_{z,w\rho}$, $V_{\Poy_{z,w\rho}}=\delta_w$. Therefore we have identified the extreme points and this prooves Theorem~\ref{thm:results:t-spec}. As a direct consequence we get

\korollar{
If $\Pr\in C(\pi^\Gf)$, then the sequence $(\xi_{B_k})_k$ satisfies a law of large numbers,
\begin{equation*}
  \frac{\xi_{B_k}}{\rho(B_k)}\to -W\log(1-z) \qquad\Pr\fsm
\end{equation*}
If $\Pr$ is extremal, then $W=w\in\R_+$ $\Pr$-a.s, and $\Pr$ is the Pólya sum process for the parameters $z$ and $w\rho$.
}


\subsection{Total particle number ensemble \label{subsect:proofs:tpe}}

$\Hsig_{B}$ gives information about the total number of particles inside the bounded set $B$. The kernel given by conditioning the Pólya sum process on the number of points inside $B$,
\begin{equation*}
  \pi^{\Hf}_{B}(\mu,\phi) \coloneqq \Poy_{z,\rho}\bigl( \phi|\Hsig_{B} \bigr)(\mu) = \Poy_{z,\rho}\bigl( \phi(\,\cdot\, +\mu_{B^{c}}) | \zeta_{B}=\zeta_{B}\mu\bigr),
\end{equation*}
is again a local specification and we get immediately for $\F_B$-measurable, non-negative $\phi$ on $\{\zeta_B=m\}$
\begin{equation} \label{eq:proofs:b-spec}
  \pi^{\Hf}_{B}(\,\cdot\,,\phi) = \frac{1}{\rho(B)^{[m]}} \begin{multlined}[t] 
           \int_{B^{m}} \phi(\delta_{x_{1}}+\ldots+\delta_{x_{m}}) 
              \bigl(\rho+\delta_{x_{1}}+\ldots+\delta_{x_{m-1}}\bigr)(\d x_{m})\times \\
           \times \bigl(\rho+\delta_{x_{1}}\bigr)(\d x_{2}) \rho(\d x_{1}).
  \end{multlined}
\end{equation} 
The first step is to disintegrate the rhs of~\eqref{eq:proofs:b-spec}. For
$\mu\in\Mpm$ and $j\in\N$ let $\gamma^j_B(\mu)$ be the number of sites of $\mu$ in $B\in\Bbd$
which are occupied by points with multiplicity $j$ and define 
\begin{equation*}
   \gamma_B:\MpmX\to\Mpm(\N),\qquad \mu\mapsto\sum_{j\geq 1}\gamma^j_B(\mu)\delta_j,
\end{equation*}
then if $\id:\N\to\N$ is the identity on $\N$, we have
\begin{equation*}
  \{\zeta_B=m\}=\{\gamma_B(\id)=m\}.
\end{equation*}
Therefore we get on $\{\zeta_B=m\}$ for $\Fsig_B$-measurable, non-negative $\phi$
\begin{equation*}
  \pi^{\Hf}_{B}(\,\cdot\,,\phi) = \sum_{\gamma\in\Mpm_{(m)}(\N)}
    \Poy_{z,\rho}(\phi|\gamma_B=\gamma)\Poy_{z,\rho}(\gamma_B=\gamma|\zeta_B=m),
\end{equation*}
i.e. for $f\in F_b$ with $\supp f\subseteq B\in\Bbd$ on $\{\zeta_B=m\}$ by lemma~\ref{thm:results:iteratedkernel}
\begin{align*}
  \pi^{\Hf}_{B}(\,\cdot\,,\e^{-\zeta_f}) 
     =\frac{1}{\rho(B)^{[m]}}\sum_{\gamma\in\Mpm_{(m)}(\N)}
       \frac{m!\rho(B)^{\gamma(\N)}}{c(\gamma)} \prod_{j\geq 1}
       \left[1-\frac{1}{\rho(B)}\int 1-\e^{-jf}\d\rho\right]^{\gamma(j)}.
\end{align*}
The product is in fact finite and we have a disintegration of $\pi^{\Hf}_{B}$. 
This reformulation and the fact that $(\gamma_{B_k})_{k\in\N}$ as a random
measure on $\N$ satisfies a large deviation principle under $\Poy_{z,\rho}$
allows to identify the limit of the mixing measure as the minimiser of a
variational problem.

Let $U_{k}$ be the total number of particles in $B_{k}$ normalised by its volume $\rho(B_k)$
\begin{equation*}
  U_{k}\mu \coloneqq \frac{\zeta_{B_{k}}\mu}{\rho(B_{k})},
\end{equation*}
and set furthermore
\begin{equation*}
  \Mpm_\rho=\left\{\mu\in\MpmX: \lim_{k\to\infty} \frac{\zeta_{B_{k}}\mu}{\rho(B_{k})}\text{ exists}\right\},
\end{equation*}
and $U$ on $\Mpm_\rho$ the corresponding limit and $0$ otherwise.

\prop{ \label{thm:proofs:b-lim}
Let $\rho\in\MX$ be a diffuse and infinite measure. For any $\Pr\in C(\pi^{\Hf})$, $\Pr(\M_\rho)=1$, and for $\phi\in L^{1}(\Pr)$
\begin{equation*}
  \Pr\bigl(\phi|\Htail\bigr)(\mu)=\Poy_{Z(\mu),\rho}(\phi) \qquad \Pr\text{-a.e. }\mu
\end{equation*}
with $Z$ being the solution of the equation
\begin{equation} \label{eq:proofs:b-spec:eqZ}
  \frac{Z}{1-Z}=U.
\end{equation}
Particularly $Z\in [0,1)$.
}

\begin{proof}
Let
\begin{equation*}
  \Mpm_{\Hf}=\left\{ \mu\in\Mpm: \lim_{k\to\infty}\pi^{\Hf}_{B_k}(\mu,\,\cdot\,) \text{ exists}\right\}.
\end{equation*} 
Since on $\{\zeta_{B_k}=m\}$,
\begin{align*}
  \begin{multlined}[t]
    \pi^{\Hf}_{B_k}(\,\cdot\,,\e^{-\zeta_f}) = \sum_{\gamma\in\Mpm_{(m)}(\N)}
    \Poy_{z,\rho}(\gamma_{B_k}=\gamma|\zeta_{B_k}=m) \times\\
    \times \prod_{j\geq 1}
    \left(1-\frac{1}{\rho(B_k)}\int
    1-\e^{-jf}\d\rho \right)^{\rho(B_k)\frac{\gamma(j)}{\rho(B_k)}},
  \end{multlined}
\end{align*}
if $\mu\in\Mpm_{\rho}$, then as shown in
proposition~\ref{subsect:proofs:ldp:t-ldp} the mixing
measure converges weakly to $\delta_{\bar{\kappa}}$ with $\bar{\kappa}=\sum_{j}\frac{Z(\mu)^j}{j}\delta_{j}$ and $Z(\mu)$ being the solution of~\eqref{eq:proofs:b-spec:eqZ}, hence
the product converges and
therefore $\Mpm_{\rho}\subseteq\Mpm_{\Hf}$.

For the reverse inclusion note that the non-convergence of $\frac{\mu(B_k)}{|B_k|}$ to some finite limit contradicts the weak convergence of $\pi^{\Hf}_{B_k}(\mu,\,\cdot\,)$.
\end{proof}

By a to the previous section analogue argumentation we get that for every $\Pr\in C(\pi^\Hf)$,
\begin{equation*}
  \Pr = \int\Poy_{z,\rho}V_\Pr(\d z)
\end{equation*}
with $V_\Pr$ here being the distribution of $Z$ under $\Pr$, a probability measure on $[0,1)$. Moreover,

\korollar{
If $\Pr\in C(\pi^\Hf)$, then the sequence $(\zeta_{B_k})_k$ satisfies a law of large numbers,
\begin{equation*}
  \frac{\zeta_{B_k}}{\rho(B_k)}\to U\qquad\Pr\fsm
\end{equation*} 
If $\Pr$ is extremal, then for some constant $u\in\R_+$, $U=u$ $\Pr$-a.s. and $\Pr$ is the Pólya sum process for the parameters $z=\tfrac{u}{1+u}$ and $\rho$.
}


\subsection{The particle-and-sites ensemble \label{subsect:proofs:tpose}}

In this last discussed ensemble the information about $\xi_B$ and $\zeta_B$ are
combined. The kernel
\begin{equation*}
  \pi^{\Ef}_{B}(\mu,\phi) \coloneqq \Poy_{z,\rho}\bigl( \phi|\Esig_{B} \bigr)(\mu) = \Poy_{z,\rho}\bigl( \phi(\,\cdot\, +\mu_{B^{c}}) | \zeta_{B}=\zeta_{B}\mu, \xi_{B}=\xi_{B}\mu \bigr),
\end{equation*}
is again a local specification. By lemma~\ref{thm:results:iteratedkernel}
we again get an explicit representation for $\F_B$-measurable $\phi$ on
$\{\zeta_B=m\}\cap\{\xi_B=k\}$%
\begin{align} \label{eq:proofs:tb-spec}
  \pi^{\Ef}_{B}(\,\cdot\,,\phi) = \frac{1}{{m \brack k}\rho(B)^{k}} \begin{multlined}[t] 
           \sum_{\stackrel{\gamma\in\Mpm_{(m)}(\N)}{\gamma(\N)=k}}
	   \frac{m!}{c(\gamma)} \int_{B^k} 
            \phi\bigl(i_1\delta_{x_1}+\ldots+i_k\delta_{x_k} \bigr)\rho(\d
	    x_1)\cdots\rho(\d x_k)
  \end{multlined}
\end{align}
where ${m \brack k}$ is the number of permutations of $m$ elements with exactly
$k$ cycles and $i_1\leq\ldots\leq i_k$ is the unique increasing sequence with for
each $j\in\N$ $\gamma(j)$ of the $i_n$'s equal to $j$.

From the previous section we keep the 
random variables $U_{k}$ and furthermore let $V_k$ be the number of occupied sites in $\mu$ normalised by the volume of $B_k$,
\begin{equation*}
V_{k}\mu \coloneqq \frac{\xi_{B_{k}}\mu}{\rho(B_{k})}.
\end{equation*}
Let furthermore $\Mpm_\rho$ be the set of those configurations $\mu$, where both ratios converge,
\begin{equation*}
  \Mpm_\rho=\left\{\mu\in\MpmX: \lim_{k\to\infty} \frac{\xi_{B_{k}}\mu}{\rho(B_{k})}\text{ and } \lim_{k\to\infty} \frac{\zeta_{B_{k}}\mu}{\rho(B_{k})} \text{ exist}\right\}
\end{equation*}
and denote by $U$ and $V$, respectively, on $\Mpm_\rho$ the corresponding limits and $0$ otherwise.

\prop{ \label{thm:proofs:tb-lim}
For any $\Pr\in C(\pi^\Ef)$, $\Pr(\M_\rho)=1$, and for $\phi\in L^{1}(\Pr)$,
\begin{equation}
  \Pr\bigl(\phi|\Etail\bigr)(\mu)=\Poy_{Z(\mu),W(\mu)\cdot\rho}(\phi) \qquad \Pr\text{-a.e. }\mu
\end{equation}
where $Z$ and $W$ are determined by
\begin{align*}
  W\frac{Z}{1-Z}=U \qquad
  -W\log (1-Z)=V.
\end{align*}
}
Note that if for some configuration $\mu$, either $U(\mu)=0$ or $V(\mu)=0$, then both vanish simultaneously. In this case put unambigously $Z(\mu)=W(\mu)=0$.
\begin{proof}
Repeat the arguments of the proof of proposition~\ref{thm:proofs:b-lim} to
obtain
\begin{equation*}
  Q_\mu=\Poy_{Z(\mu),W(\mu)\cdot\rho}.\qedhere
\end{equation*}
\end{proof}

For the particle-and-sites ensemble, every $\Pr\in C(\pi^\Hf)$ has a representation
\begin{equation*}
  \Pr = \int\Poy_{z,w\rho}V_\Pr(\d z,\d w)
\end{equation*}
and $V_\Pr$ is the distribution of $(Z,W)$ under $\Pr$. In the same way

\korollar{
If $\Pr\in C(\pi^\Ef)$, then the sequence $(\xi_{B_k},\zeta_{B_k})_k$ satisfies a law of large numbers,
\begin{equation*}
 \left(\frac{\xi_{B_k}}{\rho(B_k)},\frac{\xi_{B_k}}{\rho(B_k)}\right)\to (V,U)\qquad\Pr\fsm
\end{equation*}
If $\Pr$ is extremal, then $(V,U)=(v,u)$ $\Pr$-a.s. for some constants $v,u\in\R_+$, and $\Pr$ is the Pólya sum process for the parameters $z$ and $w\rho$ with $z$ and $w$ being the solution of
\begin{equation*}
  w\tfrac{z}{1-z}=u,\qquad -w\log(1-z)=v.
\end{equation*}
}






\subsection{Large Deviations \label{subsect:proofs:ldp}}

This last part contains the proofs of the weak convergence of $\pi^{\Hf}_{B_k}(\mu,\,\cdot\,)$ and $\pi^{\Ef}_{B_k}(\mu,\,\cdot\,)$ as $k\to\infty$. Their weak limits are determined by the minimiser of the variational problems below. 

\prop{\label{subsect:proofs:ldp:t-ldp}
For $\mu\in\MpmX$ such that $u=\lim_{k\to\infty}\frac{\zeta_{B_k}\mu}{|B_k|}\in\R$ exists, if $z_u$ is a solution of
\begin{equation*}
  \sum_{j\geq 1} z^j=u,
\end{equation*}
and $\bar{\kappa}=\sum_{j\geq 1} \frac{z_u^j}{j}\delta_j$, then
\begin{equation*}
  \Poy_{z,\rho}\biggl(\frac{\gamma_{B_k}}{|B_k|}\in\,\cdot\,\Bigl|\zeta_{B_k}=\mu(B_k)\biggr)\to \delta_{\bar{\kappa}}.
\end{equation*}
}

\prop{\label{subsect:proofs:ldp:bt-ldp}
For $\mu\in\MpmX$ such that $u=\lim_{k\to\infty}\frac{\zeta_{B_k}\mu}{|B_k|}\in\R$ and $v=\lim_{k\to\infty}\frac{\xi_{B_k}\mu}{|B_k|}\in\R$ exist, if $(z_{u,v},w_{u,v})$ is a
solution of the system
\begin{equation*}
  w\sum_{j\geq 1} z^j=u,\qquad w\sum_{j\geq 1} \frac{z^j}{j}=v,
\end{equation*}
and $\bar{\kappa}=w_{u,v}\sum_{j\geq 1} \frac{z_{u,v}^j}{j}\delta_j$, then
\begin{equation*}
  \Poy_{z,\rho}\biggl(\frac{\gamma_{B_k}}{|B_k|}\in\,\cdot\,\Bigl|\zeta_{B_k}=\zeta_{B_k}\mu, \xi_{B_k}=\xi_{B_k}\mu\biggr)\to \delta_{\bar{\kappa}}.
\end{equation*}
}
These results are a direct consequence of two following lemmas giving large deviation bounds and the solution of the corresponding minimisation problems. The remaining notation
follows these lemmas.

\lemma{\label{prop:proofs:upper-bds}
Subject to the above setup, the upper bounds are given by
\begin{equation}
  \limsup_{k\rightarrow\infty} \frac{1}{\rho(B_{k})}\log 
       \Poy_{z,\rho}\bigl(\exp(-\chi_{\{\gamma_{B_k}\in C_{u_{k},v_{k}}\}})\bigr)
     \leq -\inf_{\MN}\Bigl[ I+\chi_{D_{u,v}}\Bigr]\label{eq:proofs:bt:upbd}
\end{equation}
\begin{equation}
  \limsup_{k\rightarrow\infty} \frac{1}{\rho(B_{k})}\log 
       \Poy_{z,\rho}\bigl(\exp(-\chi_{\{\gamma_{B_k}\in C_{u_{k}}\}})\bigr)
     \leq -\inf_{\MN}\Bigl[ I+\chi_{D_{u}}\Bigr].\label{eq:polya:MDBd:b:upbd}
\end{equation}
}

\lemma{\label{prop:proofs:lower-bds}
Subject to the above setup, for each $\eps>0$ the lower bounds are given by
\begin{equation}
  \liminf_{k\rightarrow\infty} \frac{1}{\rho(B_{k})}\log 
       \Poy_{z,\rho}\bigl(\exp(-\chi_{\{\gamma_{B_k}\in C_{u_{k},v_{k}}^\eps\}})\bigr)
     \leq -\inf_{\MN}\Bigl[ I+\chi_{D_{u,v}^{\eps}}\Bigr]\label{eq:polya:MDBd:bt:lwbd}
\end{equation}
\begin{equation}
  \liminf_{k\rightarrow\infty} \frac{1}{\rho(B_{k})}\log 
       \Poy_{z,\rho}\bigl(\exp(-\chi_{\{\gamma_{B_k}\in C_{u_{k}}^\eps\}})\bigr)
     \leq -\inf_{\MN}\Bigl[ I+\chi_{ D_{u}^\eps}\Bigr].\label{eq:polya:MDBd:b:lwbd}
\end{equation}
}

The random element $\gamma_B$ in $\MpmN$ counts the number of sites in $B$ which are occupied with points of a given multiplicity.
For $B\in\Bbd$, because of the infinite divisibility of $\Poy_{z,\rho}$, the number $\gamma_B(j)$ of sites in $B$ occupied by points with multiplicity $j$ is Poisson distributed with intensity $\tfrac{z^j}{j}\rho(B)$, which follows from equation~\eqref{eq:results:infdiv-int}. Hence $\Poy_{z,\rho}\circ\gamma_B^{-1}$ is a Poisson process on $\N$ with the finite intensity measure
\begin{equation*}
  \tau_{z,B}=\rho(B)\sum_{j\geq 1}\frac{z^j}{j}\delta_j.
\end{equation*}
Recall that the measure we get by dropping the factor $\rho(B)$ is $\tau_z$. Since $\gamma_B$ is $\Poy_{z,\rho}$-a.s. an element of $\MpmN$ with finite first moment, we equip $\MN$
with the topology $\T$ generated by the at most linearly growing functions.

By Guo and Wu~\cite{lW95}, $(\gamma_{B_k})_{k\in\N}$ satisfies a large deviation 
principle under $\Poy_{z,\rho}$ with rate $\bigl(\rho(B_k)\bigr)_{k\in\N}$ and rate function $I(\,\cdot\,;\tau_z):\M(\N)\to[0,\infty]$,
\begin{equation*}
  I(\kappa;\tau_z)=\begin{cases}
    \tau_z(f\log f-f+1) & \text{ if } \kappa\ll\tau_z, f:=\frac{\d\kappa}{\d\tau_z}, f\log f-f+1\in L^{1}(\tau_z)\\
    \infty &\text{ otherwise}
  \end{cases}.
\end{equation*}
Any $\kappa$ for which $I(\kappa;\tau_z)$ is finite has a first moment.

In the situation of the ensembles in sections~\ref{subsect:proofs:tpe}~and~\ref{subsect:proofs:tpose} we fix $\mu\in\MpmX$ and obtain for the increasing sequence of bounded sets $(B_k)_{k\in\N}$ the two
sequences of real numbers given by $u_k=\frac{\zeta_{B_k}\mu}{\rho(B_k)}$ and $v_k=\frac{\xi_{B_k}\mu}{\rho(B_k)}$. In case of the total particle number ensemble we assume that 
$(u_k)_{k\in\N}$ and in case of the particle and sites ensemble we assume that both sequences converge to some finite limits $u$ and $v$, respectively. We denote by
\begin{align*}
  C_{u_k,v_k}&=\{\gamma\in\MpmN: \gamma(\id)=u_k\rho(B_k), \gamma(\N)=v_k\rho(B_k)\}\\
  C_{u_k}&=\{\gamma\in\MpmN: \gamma(\id)=u_k\rho(B_k)\},
\end{align*}
and by 
\begin{align*}
  D_{u,v} &= \Bigl\{\kappa\in\MN: \kappa(\id)=u, \kappa(\N)=v \Bigr\} \\
  D_{u} &= \Bigl\{\kappa\in\MN: \kappa(\id)=u \Bigr\},
\end{align*}
the point measures and measures on $\N$ with the corresponding fixed first moment and fixed total mass. Furthermore define
\begin{equation*}
  \chi_A(\kappa)=\begin{cases}
    0 & \text{ if } \kappa\in A\\
    +\infty & \text{ else}
  \end{cases}.
\end{equation*}

Both, $D_{u,v}$ and $D_{u}$, are $\T$-closed but not $\T$-open. For some function $f$ denote by  $f^{usc}$ and $f^{lsc}$ its upper and lower semicontinuous regularization, i.e. its lowest upper semicontinuous majorant and its largest lower semincontinuous minorant, respectively. We get

\lemman{Semicontinuous Regularisations of $\chi_{D_{u,v}}$}{ \label{lemma:proofs:bt-screg}
The upper and lower semicontinuous regularisations $\chi_{D_{u,v}}^{\text{usc}}$ and $\chi_{D_{u,v}}^{\text{lsc}}$ of $\chi_{D_{u,v}}$ with respect to $\T$ are
\begin{align}
  \chi_{D_{u,v}}^{\text{usc}}(\kappa) = \infty,\qquad
  \chi_{D_{u,v}}^{\text{lsc}}(\kappa) = \chi_{D_{u,v}}
\end{align}
}

\lemman{Semicontinuous Regularisations of $\chi_{D_{u}}$}{ \label{lemma:proofs:b-screg}
The upper and lower semicontinuous regularisations $\chi_{D_{u}}^{\text{usc}}$ and $\chi_{D_{u}}^{\text{lsc}}$ of $\chi_{D_{u}}$ with respect to $\T$ are
\begin{align}
  \chi_{D_{u}}^{\text{usc}}(\kappa) = \infty,\qquad
  \chi_{D_{u}}^{\text{lsc}}(\kappa) = \chi_{D_{u}}
\end{align}
}%
Both results are consequences of the fact that whenever a sequence of measures in $\MN$ converges with respect to $\T$, their total mass and their first moment need to converge,
too. From these two lemmas and~\cite[2.1.7]{DS00}, for each of the ensembles the upper bounds in proposition~\ref{prop:proofs:upper-bds} follow directly without replacing $D_{u,v}$ and $D_u$ by their
lower semicontinuous regularisations.

Before we study the lower bound, we solve the minimisation problems.

\lemman{Minimiser of $I+\chi_{D_{u,v}}$}{ \label{prop:proofs:bt-upbd-partition}
Let $0<v<u<\infty$ and $z_{u,v}$, $w_{u,v}$ be the solution of the system
\begin{equation} \label{eq:polya:MDBd:bt:defeq}
  w\sum_{j\geq 1} \frac{z^{j}}{j}=v, \qquad w\sum_{j\geq 1} z^{j}=u.
\end{equation}
Then the minimiser of $\inf_{\MN}\Bigl[I+\chi_{D_{u,v}}\Bigr]$ is given by
\begin{align}
  \bar{\kappa}=w_{u,v}\sum_{j\geq 1} \frac{z_{u,v}^{j}}{j}\delta_j.\label{eq:polya:bt:min-lsc}%
\end{align}
}

\begin{proof}
Let $\tilde{z}, w>0$, then
\begin{align*}
  I(\kappa;\tau_z) &-\sum_{j\geq 1} j\kappa(j)\log \tilde{z} - \sum_{j\geq 1} \kappa(j)\log w \\
      &= \sum_{j\geq 1} \kappa_{j}\biggl(\log\frac{\kappa(j)}{\tilde{z}^{j} w \tau_z(j)}-1\biggr) +\tau_z(\N),
\end{align*}
which has a unique minimiser on $\MN$, $\bar{\kappa}=\sum_{j\geq 1}w_{u,v} \frac{z_{u,v}^{j}}{j}\delta_{j}$ with $z_{u,v}$, $w_{u,v}$ being the solution of equation 
system~\eqref{eq:polya:MDBd:bt:defeq}. The uniqueness of the solution of the equation system~\eqref{eq:polya:MDBd:bt:defeq} can be seen by noting that 
\begin{align} \label{eq:polya:MDBd:cond-dom}
    f:(0,1)\times(0,\infty)&\ra \{ (s,t)\in\R^{2}_{+}: s>t \},\\
    (z,w)&\mapsto \Bigl(w\frac{z}{1-z},-w\log(1-z)\Bigr),
\end{align}
is injective.
\end{proof}

In a similar fashion the minimisation problem for the total particle ensemble is solved, and due to the missing condition on the number of occupied sites, $w$ drops out. Therefore

\lemman{Minimiser of $I+\chi_{D_{u}}$}{ \label{prop:polya:b:upbd-partition}
Let $z_{u}$ be the solution of
\begin{equation} \label{eq:polya:MDBd:b:defeq}
  \sum_{j\geq 1} z^{j}=u.
\end{equation}
Then the minimiser of $\inf_{\MN}\Bigl[I+\chi_{D_{u}}^{\text{lsc}}\Bigr]$ is given by
\begin{equation}
  \bar{\kappa}=\sum_{j\geq 1} \frac{z_{u}^{j}}{j}\delta_j.\label{eq:polya:b:min-lsc}%
\end{equation}
}

Since the upper semicontinuous regularisations of $\chi_{D_{u,v}}$ and $\chi_{D_u}$ are infinite, we cannot conclude directly. The Boltzmann principle~\cite{RZ93} is 
a way out:
$C_{u_k,v_k}$ and $C_{u_k}$ are replaced by $\eps$-blow-ups which are $\T$-open and shrink as $\eps\to 0$ to $D_{u,v}$ and $D_u$,
\begin{align*}
  D_{u,v,k}^{\eps} &\coloneqq \left\{\gamma\in\MN: \frac{\gamma(\id)}{\rho(B_k)}\in(u-\eps,u+\eps),                                                     \frac{\gamma(\N)}{\rho(B_k)}\in(v-\eps,v+\eps) \right\} \\
  D_{u,k}^{\eps} &\coloneqq \left\{\gamma\in\MN: \frac{\gamma(\id)}{\rho(B_k)}\in(u-\eps,u+\eps) \right\}.
\end{align*}
For all $\eps>0$ and $k$ large enough, $D_{u,v,k}^\eps\supseteq C_{u_k,v_k}$ and $D_{u,k}^\eps\supseteq C_{u_k}$. Since by the non-negativity of $\chi_\cdot$ the conditions
\begin{align}
  \lim_{L\rightarrow \infty} &\limsup_{k\rightarrow\infty} \frac{1}{\rho(B_{k})}\log 
           \Poy_{z,\rho}\Bigl(\exp\bigl( -\chi_{ \{\gamma_{B_k}\in D_{u,v,k}^{\eps}\} } \bigr)1_{\{\chi_{ \{\gamma_{B_k}\in D_{u,v,k}^{\eps}\} }\leq -L\}}\Bigr) = -\infty,\\
  \lim_{L\rightarrow \infty} &\limsup_{k\rightarrow\infty} \frac{1}{\rho(B_{k})}\log 
           \Poy_{z,\rho}\Bigl(\exp\bigl(-\chi_{ \{\gamma_{B_k}\in D_{u,k}^{\eps}\} }\bigr)1_{\{\chi_{ \{\gamma_{B_k}\in D_{u,k}^{\eps}\} }\leq -L\}}\Bigr) = -\infty,
\end{align}
are satisfied and we get for each $\eps>0$ by~\cite[Lemma 2.1.8]{DS00} the lower bounds in proposition~\ref{prop:proofs:lower-bds}.

\propn{Minimiser of $I+\chi_{D_{u,v}^{\eps}}$}{ \label{prop:polya:bt:lwbd-partition}
Let $0<v<u<\infty$. For sufficiently small $\eps>0$ there exists a pair $(z_{u,v,\eps},w_{u,v,\eps})\in(0,1)\times(0,\infty)$ such that the infimum of $I+\chi_{D_{u,v}^{\eps}}$ on $\MN$ is attained at
\begin{align}
  \bar{\kappa}_{\eps}=w_{u,v,\eps}\sum_{j\geq 1}
  \frac{z_{u,v,\eps}^{j}}{j}\delta_j\label{eq:polya:bt:min-usc}.%
\end{align}
As $\eps\ra 0$, $z_{u,v,\eps}\ra z_{u,v}$, $w_{u,v,\eps}\ra w_{u,v}$ and 
\begin{equation*}
  \T\text{-}\lim_{\eps\ra 0}\bar{\kappa}_{\eps}=\bar{\kappa},
\end{equation*}
where $\bar{\kappa}$ is given by equation~\eqref{eq:polya:bt:min-lsc} and $z_{u,v}$ and $w_{u,v}$ by equation~\eqref{eq:polya:MDBd:bt:defeq}.
}

\begin{proof}
In the proof of proposition~\ref{prop:proofs:bt-upbd-partition} we showed that for fixed $u>v$, the minimiser of $I$ on $D_{u,v}$ was given by
\begin{equation*}
  \bar{\kappa}=\bar{\kappa}(z,w)= w\sum_{j\geq 1} \frac{z^{j}}{j}\delta_j
\end{equation*}
with $(z,w)=f^{-1}(u,v)$ and $f$ given in equation~\eqref{eq:polya:MDBd:cond-dom}. Since the mapping $(z,w)\mapsto\bar{\kappa}(z,w)$ is continuous wrt. $\T$, it suffices to note
that $f$ and $f^{-1}$ are continuous. Therefore we get the existence of the minimiser and as $\eps\to 0$, the family of minimisers converges to the desired limit since
$D_{u,v}=\bigcap_{\eps>0}D_{u,v}^\eps$.

\end{proof}

\propn{Minimiser of $I+\chi_{D_{u}^{\eps}}$}{ \label{prop:polya:b:lwbd-partition}
Let $0\leq u<\infty$. Then there exists $z_{u,\eps}\in[0,\infty)$ such that the infimum of $I+\chi_{D_{u}}^{\eps}$ on $\MN$ is attained at
\begin{align}
  \bar{\kappa}_{\eps}=\sum_{j\geq 1} \frac{z_{u,\eps}^{j}}{j}\delta_j\label{eq:polya:b:min-usc}.%
\end{align}
As $\eps\ra 0$, $z_{u,\eps}\ra z_{u}$ and 
\begin{equation*}
  \T\text{-}\lim_{\eps\ra 0}\bar{\kappa}_{\eps}=\bar{\kappa},
\end{equation*}
where $\bar{\kappa}$ is given by equation~\eqref{eq:polya:b:min-lsc} and $z_{u}$ by equation~\eqref{eq:polya:MDBd:b:defeq}.

}

\emph{I am very grateful for the referee's remarks which lead to several improvements, in particular for the hint to~\cite{GW82}.}


\end{document}